\documentclass[11pt]{amsart}
\usepackage{amssymb, amsthm, amscd, amsmath}
\usepackage[vcentermath]{youngtab}

\newtheorem{theorem}{Theorem}[section]
\newtheorem{lemma}[theorem]{Lemma}

\theoremstyle{definition}

\newtheorem{remark}[theorem]{Remark}

\begin{document}
\renewcommand{\baselinestretch}{1.1} 
\Yboxdim10pt

\def\image{\mathop{\mathrm{Image}}\,}
\def\rank{\mathop{\mathrm{rank}}\,}
\def\dim{\mathop{\mathrm{dim}}\,}
\def\deg{\mathop{\mathrm{deg}}\,}
\def\div{\mathop{\mathrm{div}}\,}
\def\Pf{\mathrm{Pf}}
\def\Cl{\mathop{\mathrm{Cl}}\,}
\def\Pic{\mathop{\mathrm{Pic}}\,}
\def\Proj{\mathop{\mathrm{Proj}}\,}
\def\Spec{\mathop{\mathrm{Spec}}\,}
\def\Supp{\mathop{\mathrm{Supp}}\,}
\def\Hom{\mathop{\mathrm{Hom}}\,}
\def\tor{\mathrm{Tor}}
\def\CH{\mathop{CH}}
\def\ch{\mathop{\mathrm{ch}}}
\def\td{\mathop{\mathrm{td}}}
\def\chinf{\mathop{\chi_{\infty}}}
\def\Fdot{F_{\bullet}}
\def\Gdot{G_{\bullet}}
\def\m{\mathfrak{m}}
\def\p{\mathfrak{p}}
\def\sgn{\mathop{\mathrm{sgn}}}

\title{Todd classes of affine cones of Grassmannians} 

\author{Kazuhiko Kurano \and Anurag K. Singh}

\address{
\center{
Kazuhiko Kurano \newline
Department of Mathematics \newline 
Tokyo Metropolitan University \newline
Minami-Ohsawa 1-1, Hachioji \newline
Tokyo 192-0397, Japan \newline
E-mail: {\tt kurano@comp.metro-u.ac.jp} \newline
\phantom{.} \newline
Anurag K. Singh \newline
Department of Mathematics \newline
University of Utah \newline
155 S. 1400 E. \newline
Salt Lake City, UT 84112, USA \newline
E-mail: {\tt singh@math.utah.edu}}}

\thanks
{Mathematics Subject Classification Primary 14C17; Secondary 13H15, 14C15, 
14C40. \\
The first author is supported in part by the Grants-in-Aid for Scientific
Research, The Ministry of Education, Science and Culture, Japan. The second
author is supported in part by the National Science Foundation under Grant No.
DMS 0070268}

\begin{abstract} A ring $A$ which is a homomorphic image of a regular local
ring $S$ is said to be a {\it Roberts ring}\/ if $\tau_{A/S}([A]) = [\Spec
A]_{\dim A}$, where $\tau_{A/S}$ is the Riemann-Roch map for $\Spec A$. Such
rings satisfy a vanishing theorem for intersection multiplicities, as was
proved by P.~Roberts. It is known that complete intersections are Roberts
rings, and the first author showed that a determinantal ring is a Roberts ring
precisely if it is a complete intersection. Let $A_d(n)$ denote the affine cone
of the Grassmann variety $G_d(n)$ under the Pl\"ucker embedding. In this paper,
we determine precisely when $A_d(n)$ is a Roberts ring. \end{abstract}

\maketitle

\section{Introduction} 

Let $M$ and $N$ be two modules over a local ring $A$ such that $M\otimes_A
N$ has finite length and $M$ has finite projective dimension. In \cite{Serre} 
Serre defined the intersection multiplicity of $M$ and $N$ as 
$$
\chi(M,N) = \sum_{i=0}^{\dim A}(-1)^i \ell(\tor_i^A(M,N))
$$
where $\ell(-)$ is the length function. If $A$ is a regular local ring, Serre 
showed that the condition $\ell (M\otimes_A N) < \infty$ implies 
$\dim M + \dim N \le \dim A$, and posed the following conjectures:

\begin{tabbing}
\hspace{0.8in} \= \hspace{3in} \kill
{\bf Vanishing:} \> If $\dim M + \dim N < \dim A$, then $\chi(M,N) = 0$. \\
\\
{\bf Positivity:} \> If $\dim M + \dim N = \dim A$, then $\chi(M,N) > 0$. \\
\end{tabbing}
Serre settled these conjectures affirmatively for regular local rings $A$ which
are equicharacteristic or unramified of mixed characteristic. The positivity
conjecture remains open, though Gabber recently proved that $\chi(M,N) \ge 0$,
see \cite{Be}. The vanishing conjecture was proved by Roberts in 
\cite{RobVanish} and independently by Gillet and Soul\'e in \cite{GS1}. The
theorem of Roberts is as follows:

\begin{theorem}[Roberts]
Let $A$ be a homomorphic image of a regular local ring $S$ such that 
$\tau_{A/S}([A]) = [\Spec A]_{\dim A}$ (e.g., suppose that $A$ is a regular
local ring or a complete intersection). If $M$ and $N$ are finitely generated
$A$-modules of finite projective dimension such that $\ell(M \otimes_A N) 
< \infty$ and $\dim M + \dim N < \dim A$, then $\chi(M,N) = 0$.
\end{theorem}

It draws attention to rings $A$ for which $\tau_{A/S}([A])=[\Spec A]_{\dim A}$,
and such rings were named {\it Roberts rings}\/ by the first author and studied
in \cite{RobertsRings}. Complete intersections, for example, are Roberts rings.
Our primary goal in this paper is to determine when affine cones over Grassmann
varieties are Roberts rings, and we prove the result:

\begin{theorem}\label{main}
For $1 \le d \le n-1$, let $A_d(n)$ denote the affine cone of the Grassmann
variety $G_d(n)$ under the Pl\"ucker embedding. Then $A_d(n)$ is a Roberts ring
if and only if one of the following conditions is satisfied:
\begin{enumerate}
\item $d=1$;
\item $d=n-1$;
\item $d=2$ and $n=4$;
\item $d=3$ and $n=6$.
\end{enumerate}
\end{theorem}

Theorem~\ref{main} gives plenty of examples of Gorenstein factorial rings which
are not Roberts rings. The first examples of Gorenstein rings which are not
Roberts rings were discovered by the first author in \cite{KuranoTohoku} where
he computed the Todd classes of certain determinantal rings. A few years later
the second author, in collaboration with C.~Miller, \cite{MS}, found a
Gorenstein ring which is not a {\it numerically Roberts ring}\/ in the sense of
\cite{numRob}. (A local ring is a numerically Roberts ring if and only if  the
Dutta multiplicity  coincides with the Euler characteristic for any bounded
free complex with homology of finite length. We remark that a Roberts ring is a
numerically Roberts ring but the converse is not true, see \cite{numRob}.) 
Recently Roberts and Srinivas, using a localization sequence in $K$-theory,
established the existence of large families of Gorenstein rings which are not
numerically Roberts rings (\cite{RobertsSrinivas}) and, applying their methods,
we know that $A_d(n)$ is a Roberts ring if and only if it is a numerically
Roberts ring. Therefore Theorem~\ref{main} gives many examples of Gorenstein
factorial rings which are not numerically Roberts rings.

As a corollary of our results, we show in \S \ref{pfaff} that rings defined by
Pfaffian ideals are Roberts rings if and only if they are complete
intersections. The first author had earlier established that determinantal
rings are Roberts rings if and only if they are complete intersections,
\cite[6.2]{RobertsRings}. What is most curious in our theorem above, is that
the ring $A_3(6)$ is not a complete intersection, yet it is a Roberts ring.

\section{Background}

We first review some notation and results from \cite{Fu, KuranoTohoku,
RobertsRings} that we use later in our work.

\subsection*{Roberts rings}

Let $A$ be a homomorphic image of a regular local ring $S$, and $d=\dim A$. 
The {\it Chow group}\/ of $A$ is 
$$
A_*(A) = \bigoplus_{i=0}^{d} A_i(A),
$$
where $A_i(A)$ is the free abelian group generated by cycles of the form $[A/P]$
for $P \in \Spec A$ with $\dim A/P = i$, considered modulo rational equivalence.
Let $G_0(A)$ be the Grothendieck group of finitely generated $A$-modules. For 
an abelian group $M$, we use $M_{\mathbb Q}$ to denote the tensor product 
$M \otimes_{\mathbb Z} {\mathbb Q}$. With this notation, consider the 
{\it Riemann-Roch map}\/ as in \cite[Chapter 18]{Fu},
$$
\tau_{A/S}: G_0(A)_{\mathbb Q} \to A_*(A)_{\mathbb Q}.
$$
This is an isomorphism of ${\mathbb Q}$-vector spaces and it is known that
under mild hypotheses (e.g., if $A$ is complete, or is essentially of finite
type over a field or over ${\mathbb Z}$) it does not depend on the choice of
the regular local ring $S$, see \cite{RobertsRings, RobBook}. When $\tau_{A/S}$
does not depend on the choice of $S$, we denote it simply by $\tau_A$. Let $[A]$
denote the class of the ring $A$ in $G_0(A)_{\mathbb Q}$. Then $A$ is said to 
be a {\it Roberts ring}\/ if 
$$
\tau_{A/S}([A]) \in A_d(A)_{\mathbb Q}
$$
for some choice of $S$. In other words, if we write
$$
\tau_{A/S}([A]) = \tau_d + \tau_{d-1} + \cdots + \tau_0 
\quad \text{for} \quad \tau_i \in A_i(A)_{\mathbb Q},
$$
then $A$ is a Roberts ring if and only if $\tau_{d-1} = \cdots = \tau_0 = 0$ 
for some choice of $S$. We summarize some properties of Roberts rings; see 
\cite{RobertsRings} for the proofs.

\begin{theorem}\label{basic}
Consider a ring $(A,\m)$ which is a homomorphic image of a regular local ring.
\begin{enumerate}

\item If $(A,\m)$ is a Roberts ring, so are the local rings $A_{\p}$ for $\p
\in \Spec A$, and the $\m$-adic completion $\hat{A}$. 

\item Let $x\in\m$ be a nonzerodivisor. If $A$ is a Roberts ring, so is $A/xA$.

\item If $A$ is a subring of a regular local ring $T$ such that $T$ is a
module-finite extension of $A$, then $A$ is a Roberts ring.

\item If $A$ is a normal domain with a Noether normalization $S \subseteq A$
such that the extension $A/S$ is generically Galois, then $A$ is a Roberts ring.
\end{enumerate}
\end{theorem}

We next record some facts about $\tau_{A/S}$.

\begin{theorem} 
Let $A$ be a local ring of dimension $d$ which is a homomorphic image of a 
regular local ring $S$. Let
$$
\tau_{A/S}([A]) = \tau_d + \tau_{d-1} + \cdots + \tau_0 
\quad \text{for} \quad \tau_i \in A_i(A)_{\mathbb Q}.
$$

\begin{enumerate}
\item $\tau_d \neq 0$.

\item If $A$ is a complete intersection, then $\tau_i = 0$ for all $i < d$.

\item If $A$ is a Cohen-Macaulay ring with canonical module $\omega_A$, then 
$$
\tau_{A/S}([\omega_A])=\tau_d -\tau_{d-1} +\tau_{d-2} -\cdots + (-1)^d \tau_0.
$$

\item If $A$ is a Gorenstein ring, then $\tau_{d-i} = 0$ for odd integers $i$.

\item If $A$ is a normal Roberts ring, then it is ${\mathbb Q}$-Gorenstein. A
normal domain $A$ of dimension two is a Roberts ring if and only if $A$ is a
${\mathbb Q}$-Gorenstein ring.
\end{enumerate}
\end{theorem}

\subsection*{Affine cones of smooth projective varieties}

Let $R = \oplus_{n \ge 0} R_n$ be a graded ring over a field $R_0=K$ which is
generated, as a $K$-algebra, by finitely many elements of degree one. Let $\m$
be the unique homogeneous maximal ideal of $R$. Assume that $X=\Proj R$ is a
smooth projective variety of dimension $t$. Let 
$$
A_*(X) = \bigoplus_{i=0}^t A_i(X) \qquad \text{and} \qquad
\CH(X) = \bigoplus_{i=0}^t{\CH}^i(X)
$$
denote the {\it Chow group}\/ and the {\it Chow ring}\/ of $X$ respectively,
where 
$$
{\CH}^i(X)=A_{t-i}(X) \qquad \text{for all} \qquad 0 \le i \le t.
$$ 
Set ${\CH}(X)_{\mathbb Q} = {\CH}(X) \otimes {\mathbb Q}$. Let 
$h=c_1({\mathcal O}_X(1)) \cap [X] \in {\CH}^1(X)_{\mathbb Q}$ be the first
Chern class of the invertible sheaf ${\mathcal O}_X(1)$. One of the main results
of \cite{KuranoTohoku} is the following theorem:

\begin{theorem}\label{sequence}
There is an exact sequence of graded modules
$$
\begin{CD}
\CH(X)_{\mathbb Q} @>h>> \CH(X)_{\mathbb Q} @>\xi>> A_*(R_{\m})_{\mathbb Q} 
@>>> 0,
\end{CD}
$$
where $\xi$ is a map satisfying $\xi([\Proj R/\p]) = [R_{\m}/\p R_{\m}]$ for 
each homogeneous prime ideal $\p$ of $R$. Under this map,
$\xi(\td (\Omega^{\vee}_X)) = \tau_{R_{\m}}([R_{\m}])$ where 
$\td (\Omega^{\vee}_X)$ is the Todd class of the tangent sheaf $\Omega^{\vee}_X$
and $\tau_{R_{\m}}:G_0(R_{\m})_{\mathbb Q}\to A_*(R_{\m})_{\mathbb Q}$ is the 
Riemann-Roch isomorphism. In particular, 
\begin{align*}
A_0(R_{\m})_{\mathbb Q}     &\cong 0, \\   
A_{i}(R_{\m})_{\mathbb Q}   &\cong 
  \frac{{\CH}^{t+1-i}(X)_{\mathbb Q}}{h{\CH}^{t-i}(X)_{\mathbb Q}} 
  \quad \text{for all} \quad 1 \le i \le t, \quad \text{and}, \\
A_{t+1}(R_{\m})_{\mathbb Q} & \cong {\CH}^0(X)_{\mathbb Q}
\cong {\mathbb Q}.
\end{align*}
\end{theorem}

\subsection*{Dutta multiplicity}

Let $A$ be a complete local ring of dimension $d$ over a perfect field of prime
characteristic $p$, and let $\Gdot$ be a bounded complex of free modules with
homology of finite length. We denote by $F^n(-)$ the $n$\,th iteration of the
Frobenius functor. The {\it Dutta multiplicity}\/ of $\Gdot$ is the limit
$$ 
\chinf(\Gdot) = \lim_{n \to \infty}\frac{\chi(F^n(\Gdot))}{p^{nd}}, 
$$
studied by S.~Dutta in \cite{FrobMult}. The Dutta multiplicity behaves, in 
many ways, better than the usual multiplicity, and Roberts used the Dutta
multiplicity in an essential way in his proof of the new intersection theorem
in mixed-characteristic, \cite{RobIntersect}. While we do not pursue it here,
the Dutta multiplicity can be defined in a characteristic-free way, as was 
accomplished by the first author in \cite{DuttaMult}. One of the motivating
reasons for the study of Roberts rings is that over these rings the Dutta
multiplicity of a complex coincides with its Euler characteristic. 

By the support of a complex $\Gdot$, denoted $\Supp(\Gdot)$, we mean the union
of the supports of its homology modules. We summarize some results from
\cite{Jean, Fu, KR, RobVanish, Ri} which illustrate the behavior of
$\chinf(\Gdot)$.

\begin{theorem}\label{Dutta}
Let $\Gdot$ and $\Fdot$ be bounded complexes of finitely generated free modules
over a complete local ring $A$ of dimension $d$ and characteristic $p>0$. 
Assume furthermore that $A$ has a perfect residue class field.
\begin{enumerate}

\item If $\Supp(\Gdot)=\{\m\}$, then $\chinf({\Gdot}^\vee)=(-1)^d\chinf(\Gdot)$.

\medskip

\item If $\Supp(\Gdot)=\{\m\}$, then $\chinf(F^n(\Gdot))= p^{nd}\chinf(\Gdot)$
for all $n \in {\mathbb N}$.

\medskip

\item If $\dim \Supp(\Gdot) + \dim \Supp(\Fdot) < d$ and $\Supp(\Gdot) \cap
\Supp(\Fdot) = \{\m\}$, then $\chinf(\Gdot \otimes_A \Fdot) = 0$.

\medskip

\item If $\dim \Supp(\Gdot) + \dim \Supp(\Fdot) \leq d$ and $\Supp(\Gdot) \cap
\Supp(\Fdot) = \{\m\}$, then $\chinf(\Gdot \otimes_A \Fdot) = (-1)^{d - \dim
\Supp(\Gdot)}\chinf(\Gdot^\vee \otimes_A \Fdot)$.

\medskip

\item If $\Gdot$ has length $d$ and $\Supp(\Gdot)=\{\m\}$ (in particular, 
$\Gdot$ is not exact), then $\chinf(\Gdot) > 0$.
\end{enumerate}
\end{theorem}

The assertions (1), (2), (4) and (5) of Theorem~\ref{Dutta} are not true in
general if the Dutta multiplicity $\chinf$ is replaced by the usual Euler
characteristic $\chi$. However if the ring $A$ is a Roberts ring, all
assertions of Theorem 2.4 are true for the Euler characteristic $\chi$ since,
in this case, $\chinf(\Gdot) = \chi(\Gdot)$ for a bounded free complex $\Gdot$
with $\Supp(\Gdot)=\{\m\}$. 

\begin{remark} Assume that $A$ is a $d$-dimensional local ring (not necessary
of positive characteristic) which is a homomorphic image of a regular local
ring. In this generality, the Dutta multiplicity of complexes with support in
$\{\m\}$ is defined in \cite{DuttaMult}.

The statements (1), (3), (4) in Theorem~\ref{Dutta} hold true in this case.
Furthermore, if we replace $F$ with the $p$\,th adams operation $\psi^p$,
\cite{GS2,KR}, statement (2) is also valid for any positive integer $p$. If $A$
contains a field, then (5) is true, see \cite{KR}. However, this is an open
problem in the case of mixed-characteristic. The positivity of the Dutta
multiplicity is deeply connected to the positivity conjecture of Serre, 
\cite[Theorem~1.2]{Test}. \end{remark}

\begin{remark} The concept of a {\it numerically Roberts ring}\/ is defined in
\cite{numRob}. It is proved there that a local ring $A$ is a numerically
Roberts ring if and only if, over the ring $A$, the Dutta multiplicity always
coincides with the Euler characteristic. Consequently, a Roberts ring is a
numerically Roberts ring, but there are many examples of numerically Roberts
rings which are not Roberts rings.

However, using a method established in \cite{RobertsSrinivas}, the affine cone
$A_d(n)$ of a Grassmann variety is a Roberts ring if and only if it is a
numerically Roberts ring. A key point here is that for a Grassmann variety
$G=G_d(n)$, we have ${\CH}(G)_{\mathbb Q} \cong
{\CH}_{\mathrm{num}}(G)_{\mathbb Q}$. \end{remark}

\section{Vector bundles}

We review definitions and basic facts on Chern characters etc., \cite[3.2]{Fu},
that we use later.

Let $E$ be a vector bundle on a scheme $X$. We use $c_t(E)$ to denote its 
{\it Chern polynomial},
$$
c_t(E) = 1+c_1(E)t+c_2(E)t^2+c_3(E)t^3+\cdots.
$$
For an exact sequence of vector bundles $0 \to E' \to E \to E'' \to 0$, the
{\it Whitney sum} formula gives 
$$
c_t(E) = c_t(E')c_t(E'').
$$
If the vector bundle $E$ has rank $r$, then $c_i(E) = 0$ for all $i > r$. If 
its Chern polynomial is factored formally as
$$
c_t(E) = \prod_{i=1}^r (1+\alpha_i t),
$$
the $\alpha_i$'s are called the {\it Chern roots} of $E$, and the Chern classes 
of $E$ are elementary symmetric functions of $\alpha_1, \dots, \alpha_r$. The 
{\it Chern character} of $E$ is 
$$
\ch(E) = \sum_{i=1}^r \exp(\alpha_i) = \sum_{i=1}^r \sum_{n \ge 0}
 \frac{\alpha_i^n}{n !}.
$$
The first few terms, as can be found in \cite[Example 3.2.3]{Fu}, are
\begin{multline*}
\ch(E) = r + c_1 + \frac{1}{2}(c_1^2-2c_2) + \frac{1}{6}(c_1^3-3c_1c_2+3c_3) \\
+ \frac{1}{24}(c_1^4-4c_1^2c_2 + 4c_1c_3 + 2c_2^2 - 4c_4) + \cdots,
\end{multline*}
where $c_i = c_i(E)$.

The Chern character of a tensor product of vector bundles is
$$
\ch(E\otimes E') = \ch(E)\ch(E'),
$$
and for an exact sequence $0 \to E' \to E \to E'' \to 0$, we have
$$
\ch(E) = \ch(E')+\ch(E'').
$$
The Chern classes of the dual bundle $E^{\vee}$ are given by
$$
c_i(E^{\vee}) = (-1)^i c_i(E).
$$
The {\it Todd class}\/ $\td(E)$ of a vector bundle $E$ with Chern roots
$\alpha_1, \dots, \alpha_r$ is
$$
\td(E) = \prod_{i=1}^r \frac{\alpha_i}{1-\exp({-\alpha_i})},
$$
and the first few terms of the expansion are
\begin{multline*}
\td(E) = 1 + \frac{1}{2}c_1 + \frac{1}{12}(c_1^2+c_2) + \frac{1}{24}(c_1c_2) \\
+ \frac{1}{720}(-c_1^4+4c_1^2c_2 + 3c_2^2 +c_1c_3-c_4) + \cdots ,
\end{multline*}
see \cite[Example 3.2.4]{Fu}.

\section{Grassmannians}

Let $X=(x_{ij})$ be an $n \times d$ matrix of indeterminates over a field $K$, 
and consider the ring $R$ generated, as a $K$-algebra, by all the $d \times d$
minors of the matrix $X$. Then $R$ is the homogeneous coordinate ring of the
Grassmann variety $G_d(n)$ of $d$-dimensional subspaces in an $n$-dimensional
vector space, i.e., $G_d(n) = \Proj R$. The relations between the minors are
quadratic, and are the well-known Pl\"ucker relations, see \cite[Chapter VII,
\S 6]{HP}.

Setting $G=G_d(n)$ for the notational convenience, we have the universal exact 
sequence
$$
0 \to S \to {\mathcal O}^n_G \to Q \to 0
$$
where $Q$ (resp.\ $S$) is the {\em universal} rank $(n-d)$ {\em quotient
bundle} (resp.\ {\em universal} rank $d$ {\em subbundle}) on $G$, see
\cite[Chapter 14.6]{Fu}. We briefly explain the construction of $Q$ and $S$.
Consider the $K[x_{ij}]$-module $T$ which is the submodule of the free module
$K[x_{ij}]^n$ generated by the columns of $X$. Let $N$ denote the set of
elements of $T$ which have entries in $R$. Then $N$ is a graded submodule of
$R^n$, and $S$ is the locally free sheaf corresponding to $N$. Similarly, the
locally free sheaf $Q$ corresponds to the graded $R$-module $R^n/N$, see
\cite[Chapter 10.2]{RobBook}. 

By \cite[B.5.8]{Fu}, we have
$$
\Omega^{\vee}_G = \Hom(S,Q) = S^{\vee} \otimes Q. 
$$
From the universal exact sequence we also get 
$$
\wedge^n {\mathcal O}^n_G \cong {\mathcal O}_G \cong \wedge^d S \otimes 
\wedge^{n-d}Q,
$$
and so $\wedge^{n-d}Q \cong (\wedge^d S)^{\vee} \cong \wedge^d (S^{\vee})$.
Here, $\wedge^{n-d}Q$ is the very ample invertible sheaf corresponding to the
Pl\"ucker embedding $G = G_d(n) \hookrightarrow {\mathbb P}^{\binom{n}{d}-1}.$
By \cite[Remark 3.2.3 (c)]{Fu}, $c_1(\wedge^{n-d}Q)=c_1(Q)$ and setting 
$h = c_1(Q)$, Theorem~\ref{sequence} gives us the exact sequence
$$
\begin{CD}
\CH(G)_{\mathbb Q} @>h>> \CH(G)_{\mathbb Q} @>\xi>> A_*(A)_{\mathbb Q} @>>> 0
\end{CD}
$$
such that
$$
\xi(\td(\Omega^{\vee}_G)) = \tau_A([A]) \in A_*(A)_{\mathbb Q},
$$
where $A = R_{\m}$. Let $t=d(n-d)$, which is the dimension of the projective
variety $G$. Then $\dim A = t+1$, and suppose
$$
\tau_A([A]) = \tau_{t+1} + \tau_t + \tau_{t-1} + \dots + \tau_0,
$$
where $\tau_i \in A_i(A)_{\mathbb Q}$ for each $i$. Here, since $A$ is
essentially of finite type over a field, the Riemann-Roch map $\tau_{A/S}$ is
independent of the choice of a regular local ring $S$. Comparing terms with
the expansion of $\td(\Omega^{\vee}_G)$, we see that 
\begin{align*} 
\tau_{t+1} &= 1, \\
\tau_{t}   &= \frac{1}{2}c_1(\Omega^{\vee}_G) \mod h {\CH}^0(G)_{\mathbb Q}, \\
\tau_{t-1} &= \frac{1}{12}(c_1(\Omega^{\vee}_G)^2+c_2(\Omega^{\vee}_G))
    \mod h {\CH}^1(G)_{\mathbb Q}, \\
\tau_{t-2} &= \frac{1}{24}c_1(\Omega^{\vee}_G)c_2(\Omega^{\vee}_G)
    \mod h {\CH}^2(G)_{\mathbb Q}, \\
\tau_{t-3} &= \frac{1}{720}(-c_1(\Omega^{\vee}_G)^4 
  + 4c_1(\Omega^{\vee}_G)^2c_2(\Omega^{\vee}_G) + 3c_2(\Omega^{\vee}_G)^2 \\
 & \qquad\qquad\qquad + c_1(\Omega^{\vee}_G)c_3(\Omega^{\vee}_G)
  -c_4(\Omega^{\vee}_G)) \mod h {\CH}^3(G)_{\mathbb Q}, \quad \text{etc.}
\end{align*}
Recall that $A$ is a Roberts ring if and only if $\tau_i = 0$ for all $i \le t$,
and we will prove Theorem~\ref{main} essentially by establishing the vanishing
or nonvanishing of $\tau_i$'s.

\bigskip

\noindent {\it Proof of Theorem~\ref{main}.} If $d=1$ or $d=n-1$, the affine
cone $A_d(n)$ is a regular local ring, and therefore is a Roberts ring.
Consequently we may assume $2 \le d \le n-2$. In the case $d=2$ and $n=4$, it
is easily seen that there is exactly one Pl\"ucker relation, and so $A_2(4)$ is
a hypersurface, hence a Roberts ring.

In general, the Whitney sum formula, applied to the universal exact sequence
$$
0 \to S \to {\mathcal O}^n_G \to Q \to 0, \text{ \quad gives \quad}
c_t(S)c_t(Q) = c_t({\mathcal O}^n_G) = (c_t({\mathcal O}_G))^n = 1,
$$
which says
\begin{multline*}
(1+c_1(S)t + c_2(S)t^2 + c_3(S)t^3 + c_4(S)t^4 + \cdots) \\
\times (1+c_1(Q)t + c_2(Q)t^2 + c_3(Q)t^3 + c_4(Q)t^4 + \cdots) = 1.
\end{multline*}
Comparing the coefficients, we obtain
\begin{align*}
&c_2(S) +c_1(S)c_1(Q) + c_2(Q) = 0, \text{ \quad and \quad} \\
&c_4(S) +c_3(S)c_1(Q) + c_2(S)c_2(Q)+ c_1(S)c_3(Q) + c_4(Q)= 0 .
\end{align*}

By Lemma~\ref{young} (1) below, the graded component of 
$\CH(G)_{\mathbb Q}/ h\CH(G)_{\mathbb Q}$ in degree one is
$$
{\CH}^1(G)_{\mathbb Q} / h {\CH}^0(G)_{\mathbb Q}=0.
$$
In particular, $c_1(E) \equiv 0 \mod h \CH(G)_{\mathbb Q}$ for any vector 
bundle $E$ on $G$. Hence we have $c_2(S) \equiv -c_2(Q)$ and $c_4(S) \equiv
c_2(Q)^2 - c_4(Q)$, which will be used later. Furthermore, the expansion of
$\ch(E)$ is simplified as
\begin{align*}
\ch(E) & \equiv \rank E + \frac{1}{2}(-2c_2(E)) + \frac{1}{6}(3c_3(E)) 
+ \frac{1}{24}(2c_2(E)^2 - 4c_4(E)) + \cdots \\
& \equiv \rank E -c_2(E) + \frac{1}{2}c_3(E) + \frac{1}{12}(c_2(E)^2 - 2c_4(E)) 
+ \cdots .
\end{align*}
Recall that $\ch(\Omega^{\vee}_G) = \ch(S^{\vee}) \ch(Q)$ since 
$\Omega^{\vee}_G = S^{\vee}\otimes Q$. This equation gives us
\begin{align*}
d(n-d) & - c_2(\Omega^{\vee}_G) + \frac{1}{2}c_3(\Omega^{\vee}_G) 
+ \frac{1}{12}(c_2(\Omega^{\vee}_G)^2 - 2c_4(\Omega^{\vee}_G)) +\cdots \\
& \equiv \left[d - c_2(S^{\vee}) + \frac{1}{2}c_3(S^{\vee}) 
+ \frac{1}{12}(c_2(S^{\vee})^2 - 2c_4(S^{\vee})) + \cdots \right] \\
& \quad \times \left[n-d - c_2(Q) + \frac{1}{2}c_3(Q) 
+ \frac{1}{12}(c_2(Q)^2 - 2c_4(Q)) + \cdots \right] . \quad (*)
\end{align*}
Comparing the components of degree two in equation~$(*)$, we see that
$$
c_2(\Omega^{\vee}_G) \equiv dc_2(Q) + (n-d)c_2(S^{\vee}) 
 \equiv dc_2(Q) + (n-d)c_2(S).
$$
Since $c_2(S)\equiv -c_2(Q)$, we have $c_2(\Omega^{\vee}_G) \equiv(2d-n)c_2(Q)$.
Consequently
\begin{align*}
\tau_{t-1} & = \frac{1}{12}(c_1(\Omega^{\vee}_G)^2+c_2(\Omega^{\vee}_G))
 \mod h {\CH}^1(G)_{\mathbb Q} \\
           & = \frac{1}{12}(2d-n)c_2(Q) \mod h {\CH}^1(G)_{\mathbb Q}.
\end{align*}
We need the following lemma to complete the proof of Theorem~\ref{main}:

\begin{lemma}\label{young}
Let $G$ denote the Grassmann manifold $G_d(n)$ where $2 \le d \le n-2$. Then
\begin{enumerate}
\item ${\CH}^1(G)_{\mathbb Q} = h{\CH}^0(G)_{\mathbb Q}$.

\item $c_2(Q) \notin h {\CH}^1(G)_{\mathbb Q}$.

\item If $d \ge 4$ and $n=2d$, then $c_2(Q)^2 \notin h {\CH}^3(G)_{\mathbb Q}$.

\item If $d = 3$ and $n=6$, then 
$h{\CH}^i(G)_{\mathbb Q} = {\CH}^{i+1}(G)_{\mathbb Q}$ for $3 \le i \le 8$.
\end{enumerate}\end{lemma}

\bigskip

We first complete the proof of Theorem~\ref{main} using this lemma. Recall that
we may assume $2 \le d \le n-2$. Since 
$$
\tau_{t-1} = \frac{1}{12}(2d-n)c_2(Q) \mod h {\CH}^1(G)_{\mathbb Q},
$$
Lemma~\ref{young} (2) implies that $\tau_{t-1}$ is nonzero if $n \neq 2d$.
Consequently $A_d(n)$ is not a Roberts ring in this case.

We next assume that $n=2d$. Since $c_2(\Omega^{\vee}_G) \equiv(2d-n)c_2(Q)$,
we have $c_2(\Omega^{\vee}_G) \equiv 0$. Comparing the components of degree 
four in equation~$(*)$, we get
\begin{align*}
-\frac{1}{6}c_4(\Omega^{\vee}_G) & \equiv \frac{d}{12}(c_2(Q)^2 - 2c_4(Q)) 
+c_2(S^{\vee})c_2(Q) + \frac{d}{12}(c_2(S^{\vee})^2 - 2c_4(S^{\vee})) \\
& \equiv c_2(S)c_2(Q) + \frac{d}{12}(c_2(Q)^2 - 2c_4(Q) + c_2(S)^2 - 2c_4(S)).
\end{align*}
Since $c_2(S) \equiv -c_2(Q)$ and $c_4(S) \equiv c_2(Q)^2 - c_4(Q)$, we have 
$$
-\frac{1}{6}c_4(\Omega^{\vee}_G) \equiv -c_2(Q)^2. 
$$
Consequently $c_4(\Omega^{\vee}_G) \equiv 6c_2(Q)^2$ and so
$$
\tau_{t-3} = - \frac{1}{120} c_2(Q)^2 \mod h {\CH}^3(G)_{\mathbb Q}.
$$

If $n=2d$ and $d \ge 4$, then $\tau_{t-3}$ is nonzero by Lemma~\ref{young} (3).
Hence $A_d(n)$ is not a Roberts ring in this case. 

Suppose that $n = 6$ and $d= 3$. Then $A_3(6)$ is a Gorenstein ring of 
dimension $10$ and so $\tau_i =0$ for odd integers $i$. Since $n=2d$, we have 
$\tau_8 = \tau_{t-1} = 0$. The equality $\tau_{9-i} = 0$ for $i \ge 3$ follows 
from Lemma~\ref{young} (4). Hence $A_3(6)$ is a Roberts ring.

\bigskip

We now record the proof of Lemma~\ref{young}.

\begin{proof}[Proof of Lemma~\ref{young}] We shall use the notation and results
of \cite[Chapter 14.5--14.7]{Fu} for Schubert cycles. The Chow ring
$\CH(G)_{\mathbb Q}$ has a basis over ${\mathbb Q}$ represented by the set of
partitions
$$
\lambda = (\lambda_1,\dots,\lambda_d) \qquad \text{where} \qquad 
n-d \ge \lambda_1 \ge \dots \ge \lambda_d \ge 0.
$$
We denote the cycle corresponding to a partition 
$\lambda = (\lambda_1,\dots,\lambda_d)$ by $\{ \lambda \}$ or 
$\{\lambda_1,\dots,\lambda_d \}$. Set $| \lambda | = \sum\lambda_i$. Then 
${\CH}^l(G)_{\mathbb Q}$ has a basis consists of the set of cycles 
$\{\lambda \}$ such that $| \lambda | = l$. For $1\le m\le n-d$, the classes 
$c_m(Q)$ are called the {\it special Schubert classes}\/ and $\sigma_m =c_m(Q)$ 
coincides with the cycle $\{m,0,\dots,0\}$. The multiplication by $\sigma_m$ is 
determined by {\it Pieri's formula}: 
$$
\{ \lambda \} \times \sigma_m = \sum \{ \mu \} 
$$
where the sum runs over $\mu$ with 
$$
n-d \ge \mu_1 \ge \lambda_1 \ge \mu_2 \ge \lambda_2 \ge \dots \ge \mu_d 
\ge \lambda_d \quad \text{and} \quad | \mu | = | \lambda | + m.
$$ 

\bigskip\noindent
$(1)$ The group ${\CH}^1(G)_{\mathbb Q}$ is a ${\mathbb Q}$-vector space of 
dimension one, whose generator, in terms of a Young diagram, is \yng(1). Since 
$h \in {\CH}^1(G)_{\mathbb Q}$ corresponds to a very ample divisor, $h$ does 
not vanish. Therefore, ${\CH}^1(G)_{\mathbb Q} = h{\CH}^0(G)_{\mathbb Q}$ is 
satisfied.

\bigskip\noindent
$(2)$ The group ${\CH}^2(G)_{\mathbb Q}$ is a ${\mathbb Q}$-vector space of 
dimension two spanned by $\displaystyle{\yng(2) \ \text{and} \ \yng(1,1)}$.
The image $h {\CH}^1(G)_{\mathbb Q}$ is the ${\mathbb Q}$-span of 
$$
\yng(1) \times \yng(1) = \yng(1,1) + \yng(2) \ , \quad \text{and so} \quad 
c_2(Q) = \yng(2) \notin h {\CH}^1(G)_{\mathbb Q}.
$$

\bigskip\noindent
$(3)$ Since ${\CH}^3(G)_{\mathbb Q}$ is spanned by $\displaystyle{\yng(3)}$~, 
$\displaystyle{\yng(2,1)}$ and $\displaystyle{\yng(1,1,1)}$~, it follows that 
$h {\CH}^3(G)_{\mathbb Q}$ is spanned by 
$\displaystyle{\yng(3) \times \yng(1) = \yng(3,1) + \yng(4)}$ ,
$$
\yng(2,1) \times \yng(1) = \yng(3,1) + \yng(2,2) + \yng(2,1,1) \quad \text{and} 
\quad \yng(1,1,1) \times \yng(1) = \yng(2,1,1) + \yng(1,1,1,1) \ .
$$
Then it is easy to see that 
$$
c_2(Q)^2 = \yng(2) \times \yng(2) = \yng(4) + \yng(3,1) + \yng(2,2)
$$
is not an element of $h{\CH}^3(G)_{\mathbb Q}$.

\bigskip\noindent
$(4)$ If $d=3$ and $n=6$, then $h {\CH}^3(G)_{\mathbb Q}$ is spanned by 
$$
\yng(3) \times \yng(1) = \yng(3,1) \ , \quad 
\yng(2,1) \times \yng(1) = \yng(3,1) + \yng(2,2) + \yng(2,1,1) \quad \text{and} 
\quad \yng(1,1,1) \times \yng(1) = \yng(2,1,1)
$$
which, in this case, generate ${\CH}^4(G)_{\mathbb Q}$. This shows 
$h{\CH}^3(G)_{\mathbb Q} = {\CH}^4(G)_{\mathbb Q}$, and the remaining cases
may be computed similarly.
\end{proof}

\begin{remark} The ring $A_3(6)$ is not a complete intersection. It is a ring 
of dimension $10$, and is the homomorphic image of a regular local ring of
dimension $20$ (which is the number of $3 \times 3$ minors of a $6 \times 3$
matrix) modulo an ideal generated minimally by $35$ Pl\"ucker relations. The
number of minimal generators may be checked using \cite[Chapter VII, \S6]{HP}
and eliminating redundant relations, or by a computer algebra package such as
{\tt Macaulay2}. \end{remark}

\section{Pfaffian ideals}\label{pfaff}

We determine next when the rings $S/\Pf_m(Y)$ defined by Pfaffian ideals are
Roberts rings. 

Let $Z = (z_{ij})$ be an $2m \times 2m$ anti-symmetric matrix, that is, $z_{ij}
= -z_{ji}$ for $1 \leq i < j \leq 2m$ and $z_{ii} = 0$ for $1 \leq i \leq 2m$.
We call
$$
\Pf(Z) = \sum_\sigma \sgn(\sigma) z_{\sigma(1)\sigma(2)}z_{\sigma(3)\sigma(4)}
\cdots z_{\sigma(2m-1)\sigma(2m)}
$$
the {\em Pfaffian} of $Z$, where the sum is taken over permutations of
$\{ 1, 2, \ldots, 2m\}$ which satisfy 
$\sigma(1) < \sigma(3) < \cdots < \sigma(2m-1)$ and
$$
\sigma(1) < \sigma(2), \quad \sigma(3) < \sigma(4), \quad \dots, \quad 
\sigma(2m-1) < \sigma(2m).
$$
It is easy to see that $\Pf(Z)^2 = \det(Z)$.

Let $m$ and $n$ be positive integers such that $2m \leq n$, and $Y = (y_{ij})$
be the $n \times n$ anti-symmetric matrix with variables $y_{ij}$ for 
$1 \leq i < j \leq n$. For a set of integers such that 
$1 \leq s_1 < \cdots < s_{2m} \leq n$, we denote by $\Pf(s_1, \ldots, s_{2m})$ 
the Pfaffian of the $2m \times 2m$ anti-symmetric matrix $(y_{s_is_j})$. Let
$K$ be a field and $S$ be the localization of the polynomial ring 
$K[y_{ij} \mid 1 \leq i < j \leq n]$ at its homogeneous maximal ideal. We 
denote by $\Pf_m(Y)$ the ideal of $S$ generated by all the elements 
$\Pf(s_1, \ldots, s_{2m})$ for $1 \leq s_1 < \cdots < s_{2m} \leq n$. Set 
$B_m(n) = S/\Pf_m(Y)$. It is well known that $B_m(n)$ is a factorial Gorenstein 
ring and that
$$
\dim B_m(n) = \dim S - (n - 2m + 1)(n -2m + 2)/2.
$$ 
With this notation we have the following theorem:

\begin{theorem}
The following conditions are equivalent:
\begin{enumerate}
\item $B_m(n)$ is a Roberts ring;
\item $B_m(n)$ is a complete intersection;
\item $n = 2m$ or $m = 1$.
\end{enumerate}
\end{theorem}

\noindent {\it Proof.} The minimal number of generators of the ideal $\Pf_m(Y)$
is $\binom{n}{2m}$, and its height is $(n-2m+1)(n-2m+2)/2 =\binom{n-2m+2}{2}$.
Using these facts, the equivalence of (2) and (3) is easily verified.

In the case $m = 2$, the ideal $\Pf_2(Y)$ is generated by the elements
$$
y_{ij}y_{kl} - y_{ik}y_{jl} + y_{il}y_{jk}, \qquad \text{for} \qquad 
1 \leq i < j < k < l \leq n. 
$$
These are precisely the Pl\"ucker relations for the Grassmann variety $G_2(n)$,
and so $B_2(n)$ coincides with $A_2(n)$. It then follows from Theorem~\ref{main}
that $B_2(n)$ is a Roberts ring if and only if $n = 4$.

Next assume that $m \geq 3$. If $n = 2m$, then $B_m(n)$ is a complete
intersection and, therefore, a Roberts ring. If $n > 2m$, then a suitable
localization of $B_m(n)$ gives a Pfaffian ring $B_{m-1}(n-2)$ over a different
base field. By induction on $m$, we may assume that $B_{m-1}(n - 2)$ is not a
Roberts ring and it follows from Theorem~\ref{basic} (1) that $B_m(n)$ is not a
Roberts ring. This completes the proof of the theorem.

\begin{remark} The ring $A_d(n)$ is a Roberts ring if and only if it is a
numerically Roberts ring. Consequently $B_2(n)$ is a Roberts ring if and only
if it is a numerically Roberts ring. However, the authors do not know whether
or not the rings $B_m(n)$ are numerically Roberts rings in the case $m \geq 3$.
\end{remark}


\begin{thebibliography}{BFM}

\bibitem[Be]{Be} P.~Berthelot, {\em Alt\'erations de vari\'et\'es alg\'ebriques
(d'apr\`es A. J. de Jong)}, S\'em. Bourbaki Vol. 1995/96, Ast\'erisque {\bf
241} (1997), Exp. No. 815, 273--311.

\bibitem[Ch]{Jean} C-Y.~J.~Chan, {\it An intersection multiplicity in terms of
$\rm{Ext}$-modules}, Proc. Amer. Math. Soc. {\bf 130} (2002), 327--336.

\bibitem[Du]{FrobMult} S.~P.~Dutta, {\em Frobenius and multiplicities}, J.
Algebra {\bf 85} (1983), 424--448. 

\bibitem[Fu]{Fu} W.~Fulton, {\em Intersection Theory}, Second edition, 
Springer-Verlag, Berlin, 1998.

\bibitem[GS1]{GS1} H.~Gillet and C.~Soul\'e, {\em K-th\'eorie et nullit\'e des
multiplicit\'es d'intersection}, C. R. Acad. Sci. Paris S\'er. I Math. {\bf 300}
(1985), 71--74.

\bibitem[GS2]{GS2} H.~Gillet and C.~Soul\'e, {\em Intersection theory using
Adams operations}, Invent. Math. {\bf 90} (1987), 243--277.

\bibitem[HP]{HP} W.~V.~D.~Hodge and D.~Pedoe, {\em Methods of Algebraic
Geometry. Vol. I}, Cambridge, at the University Press; New York, The Macmillan
Company, 1947.

\bibitem[Ku1]{DuttaMult} K.~Kurano, {\em An approach to the characteristic free
Dutta multiplicities}, J. Math. Soc. Japan {\bf 45} (1993), 369--390. 

\bibitem[Ku2]{KuranoTohoku} K.~Kurano, {\em A remark on the Riemann-Roch
formula on affine schemes associated with Noetherian local rings}, T\^ohoku
Math. J. {\bf 48} (1996), 121--138.

\bibitem[Ku3]{RobertsRings} K.~Kurano, {\em On Roberts rings}, J. Math. Soc.
Japan {\bf 53} (2001), 333--355.

\bibitem[Ku4]{Test} K.~Kurano, {\em Test modules to calculate Dutta
Multiplicities}, J. Algebra {\bf 236} (2001), 216--235.

\bibitem[Ku5]{numRob} K.~Kurano, {\em Numerical equivalence defined on a Chow
group of a Noetherian local ring}, in preparation.

\bibitem[KR]{KR} K.~Kurano and P.~Roberts, {\em Adams operations, localized
Chern characters, and positivity of Dutta multiplicity in characteristic $0$}, 
Trans. Amer. Math. Soc. {\bf 352} (2000), 3103--3116.

\bibitem[MS]{MS} C.~M.~Miller and A.~K.~Singh, {\em Intersection multiplicities
over Gorenstein rings}, Math. Ann. {\bf 317} (2000), 155--171.

\bibitem[Ro1]{RobVanish} P.~Roberts, {\em The vanishing of intersection
multiplicities of perfect complexes}, Bull. Amer. Math. Soc. (N.S.) {\bf 13}
(1985), 127--130.

\bibitem[Ro2]{RobIntersect} P.~Roberts, {\em Le th\'eor\`eme d'intersection},
C. R. Acad. Sci. Paris S\'er. I Math. {\bf 304} (1987), 177--180.

\bibitem[Ro3]{Ri} P.~Roberts, {\it Intersection theorems}, Commutative algebra,
417--436, Math.\ Sci.\ Res.\ Inst.\ Publ. {\bf 15}, Springer, New York, Berlin,
1989.

\bibitem[Ro4]{RobBook} P.~Roberts, {\em Multiplicities and Chern classes in
local algebra}, Cambridge University Press, Cambridge, 1998.

\bibitem[RS]{RobertsSrinivas} P.~Roberts and V.~Srinivas, {\em Modules of
finite length and finite projective dimension}, Invent. Math., to appear. 

\bibitem[Se]{Serre} J.-P.~Serre, {\em Alg\`ebre
locale\,$\cdot$\,Multiplicit\'es}, Lecture Notes in Mathematics, {\bf 11}
Springer-Verlag, Berlin-New York, 1965.

\end{thebibliography}
\end{document}